\documentclass[12pt]{amsart}

\usepackage{amsmath,amssymb,amsthm,
	amsfonts}

	\newtheorem{dfn}{Definition}[section]
	\newtheorem{thm}[dfn]{Theorem}
	
	\newtheorem{cor}[dfn]{Corollary}
	\newtheorem{lem}[dfn]{Lemma}
	\newtheorem{rem}[dfn]{Remark}
	
	\newtheorem{claim}[dfn]{Claim}

	\newtheorem{ack}{Acknowledgements\!\!}

	\numberwithin{equation}{section}

	\def\notin{\not\in}

	\newcommand{\dist}{\mathop{\mathit{d}} \nolimits}
	
	\newcommand{\diam}{\mathop{\mathrm{diam}} \nolimits}
	\newcommand{\ric}{\mathop{\mathit{Ric}}  \nolimits}
	
	\newcommand{\card}{\mathop{\mathrm{Card}} \nolimits}

	\newcommand{\sep}{\mathop{\mathrm{Sep}} \nolimits}

	\newcommand{\supp}{\mathop{\mathrm{Supp}}    \nolimits}
	\newcommand{\rlip}{\mathop{\mathrm{Lip}}     \nolimits}
	\newcommand{\oblip}{\mathop{\stackrel{\rlip_1}{\longrightarrow}}    \nolimits}

\begin{document}

	\title[Observable concentration of mm-spaces into spaces with doubling measures]
	{Observable concentration of mm-spaces into spaces with doubling measures}
	\author[Kei Funano]{Kei Funano}
	\address{Mathematical Institute, Tohoku University, Sendai 980-8578, JAPAN}
	\email{sa4m23@math.tohoku.ac.jp}
	\dedicatory{}
	\date{\today}

	\maketitle

	\setlength{\baselineskip}{5mm}

	\begin{abstract}The property of measure concentration is that an arbitrary $1$-Lipschitz function $f:X\to
	 \mathbb{R}$ on an mm-space $X$ is almost close to a constant function. In this
	 paper, we prove that if such a concentration phenomenon arise, then any $1$-Lipschitz map $f$ from $X$ to a space $Y$ with a doubling
	 measure also concentrates to a constant map. As a corollary, we get any $1$-Lipschitz map to a Riemannian manifold with a lower
	 Ricci curvature bounds also concentrates to a constant map.
	\end{abstract}

	\section{Introduction}
	Let $\mu_{n}$ be the volume measure on the $n$-dimensional unit sphere $\mathbb{S}^n$ in $\mathbb{R}^{n+1}$ normalized as $\mu_n (\mathbb{S}^n)=1$. In $1919$, P. L\'{e}vy proved that
	for any $1$-Lipschitz function $f:\mathbb{S}^n \to \mathbb{R}$ and any $\varepsilon>0$, the inequality
	\begin{align*}\mu_n\big(\{ x\in \mathbb{S}^n \mid |f(x)-m_f|\geq \varepsilon    \}\big)\leq 2\ e^{-(n-1){\varepsilon}^2/2}
	 \end{align*}holds, where $m_f$ is some constant determined by $f$. For any fixed $\varepsilon>0$ the right-hand side of the above
	inequality converges to $0$ as $n\to \infty$. This
	means that any $1$-Lipschitz function on $\mathbb{S}^n$ is almost closed to a constant function for suffiecient large $n$. In 1999, M. Gromov introduced the notion of the observable diameter in \cite{gromov}. Let us recall its definition. 

	\begin{dfn}\upshape
	 Let $Y$ be a metric space and $\nu_Y$ a Borel measure on $Y$ such that $m:=\nu_Y(Y)<+\infty$. 
	 We define for any $\kappa >0$
	 \begin{align*}
	  \diam (\nu_Y , m-\kappa):= \inf \{ \diam Y_0 \mid Y_0 \subseteq Y \text{ is a Borel subset such that }\nu_Y(Y_0)\geq m-\kappa\}
	  \end{align*}and call it the \emph{partial diameter} of $\nu_Y$. 
	 \end{dfn}
	An \emph{mm-space} is a triple $(X,\dist,\mu)$, where $\dist$ is a complete separable metric on
	a set $X$ and $\mu$ a Borel measure on $(X,\dist)$ with $\mu(X)<+\infty$.

	\begin{dfn}[Observable diameter]\upshape Let $(X,\dist,\mu)$ be an mm-space and $Y$ a metric space. For any $\kappa >0$ we
	 define the \emph{observable diameter} of $X$ by 
	 \begin{align*}
	  \diam (X\oblip Y, m-\kappa):=
	   \sup \{ \diam (f_{\ast}(\mu),m-\kappa) \mid f:X\to Y \text{ is an }1 \text{{\rm -Lipschitz map}}  \},
	  \end{align*}where $f_{\ast}(\mu)$ stands for the push-forward measure of $\mu$ by $f$. 
	 \end{dfn}
	The target metric space $Y$ is called the \emph{screen}. The idea of the observable diameter came from the quantum and statistical
	mechanics, that is, we think of $\mu$ as a state on a configuration space $X$ and $f$ is interpreted as an observable. Suppose that $\diam (X\oblip \mathbb{R},m-\kappa ) <\varepsilon $ for
	suffieciently small $\varepsilon, \kappa >0$. By the definition, for any $1$-Lipschitz function $f:X\to \mathbb{R}$, there exists
	a Borel subset $A_f \subseteq \mathbb{R}$ such that $\diam A_f < \varepsilon$ and $f_{\ast}(\mu)(A_f)\geq m-\kappa$. If we pick a point $m_f\in A_f$ and fix it, then we have
	\begin{align*}
	 \mu(\{ x\in X \mid |f(x) - m_f|\geq \varepsilon        \})\leq \mu\big(f^{-1}(\mathbb{R}\setminus A_f)\big)\leq \kappa.
	\end{align*}Since $\varepsilon$ and $\kappa$ are suffieciently small positive numbers, the above inequality means that any
	$1$-Lipschitz function $f$ on $X$ is almost close to the constant function $m_f$. On the basis of this fact, we define a sequence $\{  X_n   \}_{n=1}^{\infty}$ of
	mm-spaces is a \emph{L\'{e}vy
	family} if $\diam (X_n\oblip \mathbb{R},m_n-\kappa)\to 0$ as $n\to \infty$ for any $\kappa >0$, where
	$m_n$ is the total measure of the mm-space $X_n$. Gromov proved in \cite{gromov} that if a sequence $\{  X_n
	\}_{n=1}^{\infty}$ of mm-spaces is a L\'{e}vy family, then  $\diam (X_n\oblip \mathbb{R}^k,m_n -\kappa)\to 0$ as $n\to \infty$
	for any $\kappa >0$ and $k\in \mathbb{N}$. He also discussed the case that the dimension of $\mathbb{R}^k$ goes to
	$\infty$. Our paper \cite{funano} tackles this problems in the case that the screens are the real hyperbolic spaces. Gromov treated in \cite{gromov} the case that the screen $Y$ moves around all elements of a family $\mathcal{C}_0$ of compact metric spaces which is
	precompact with respect to the Gromov-Hausdorff distance. In particular, he proves the following theorem. For any $m\in
	\mathbb{N}$, $\kappa_1\in \mathbb{R}$ and $D>0$, we denote by $\mathcal{M}_{m,\kappa_1,D}$ the set of all complete Riemannian manifolds $M$ of dimension $m$ such that
	$\ric_M \geq \kappa_1$ and $\diam M \leq D$.

	\begin{thm}[Gromov \cite{gromov}]\label{aiueo}Let $\{ X_n\}_{n=1}^{\infty}$ be a L\'{e}vy family. Then we have
	 \begin{align*}
	 \sup \{\diam (X_n \oblip M,m_n -\kappa) \mid M \in  \mathcal{M}_{m,\kappa_1,D} \} \to 0
	 \end{align*}as $n\to \infty$ for any $\kappa >0$.
	 \end{thm}
	In this paper, we consider more large class of screens. We treat the case that screen $Y$ has a doubling measure. Let $Y$
	be a metric space. Given $x\in Y$ and $r>0$, we indicate by $B_Y(x,r)$ the closed ball centered at $x$ with radius $r$. For a number $R>0$ and a function $C:(0,R]\to [1,+\infty)$, we denote by $D_{C,R}$ the set of all pairs $(Y,\nu_Y)$ satisfying
	the following properties: $\nu_Y$ is a Borel measure on $Y$ such that 
	\begin{align*}
	  0 <  \nu_Y \big( B_Y(x,r)  \big)<+ \infty
	 \end{align*}and 
		\begin{align}\label{doub}
	   \nu_Y\big(B_Y (x,2r)\big) \leq  C(r) \nu_Y \big(B_Y (x,r)  \big) \tag{1.1}
		\end{align}for all $x\in Y$ and $r>0$ with $r\leq R$. A main theorem of this paper is the following:
		 \begin{thm}\label{nobunaga}Let $\{ X_n \}_{n=1}^{\infty}$ be a L\'{e}vy family and $C:(0,R]\to [1,+\infty)$ an arbitrary
		  function. Then we have
		  \begin{align*}
		   \sup \{    \diam (X_n \oblip Y,m_n-\kappa) \mid    (Y,\nu_Y)\in D_{C,R}                  \}\to 0
		   \end{align*}as $n\to \infty$ for any $\kappa >0$.
		  \end{thm}

	Theorem \ref{nobunaga} together with Bishop-Gromov volume comparison theorem proves that Theorem \ref{aiueo} holds even
	without the diameter bound $D$ for the screens.

	\section{Preliminaries}

	\subsection{Basics of doubling measures}
	Although the following lemma and corollary are somewhat standard, we prove them for the completeness of this paper.
	\begin{lem}[cf. \cite{ambro}]\label{aaaa}Suppose that $(Y,\nu_Y)\in D_{C,R}$. Then for any $r_1,r_2>0$ with $2r_1\leq 2r_2\leq R$ , there exists a number $C(r_1,r_2)$ depending only on $r_1$ and $r_2$ such that
	 \begin{align*}
	  \frac{\nu_Y\big(  B_Y(x,r_1)  \big)}{ \nu_Y \big(   B_Y(y,r_2)         \big)} \geq \frac{1}{C(r_1,r_2)^2}
	  \Big(\frac{r_1}{r_2}\Big)^{C(r_1,r_2)}
	 \end{align*}for any $x,y\in Y$ with $x\in B_Y(y,r_2)$. 
	 \begin{proof}Put $j:= \min \{  i \in \mathbb{N}\mid B_Y(y,r_2)\subseteq B_Y(x,2^i r_1)         \}$ and $\widetilde{C}(r_1,r_2):= \max\{  C(2^i r_1)
	  \mid  1\leq i\leq   \log_2 (2r_2/r_1)     \}$. Since $B_Y(y,r_2)  \nsubseteq B_Y(x,2^{j-1}r_1)$ and $B_Y(y,r_2) \subseteq B_Y(x,2r_2)$, we have $2^{j-1}r_1
	  \leq 2r_2$. Iterating (\ref{doub}) $j$ times yields
	  \begin{align*}
	   \nu_Y \big(   B_Y(y,r_2)  \big)\leq \nu_Y \big(    B_Y(x,2^j r_1)   \big) \leq \widetilde{C}(r_1,r_2)^j \nu_Y \big(  B_Y(x,r_1)        \big).
	   \end{align*}
	  As a result, we obtain 
	  \begin{align*}
	   \frac{\nu_Y\big(   B_Y(x,r_1) \big)}{\nu_Y\big( B_Y(y,r_2)  \big)}\geq \widetilde{C}(r_1,r_2)^{-j}\geq
	   \widetilde{C}(r_1,r_2)^{-\log_2 \frac{4r_2}{r_1}}= \frac{1}{\widetilde{C}(r_1,r_2)^2}
	   \Big(\frac{r_1}{r_2} \Big)^{\frac{\widetilde{C}(r_1,r_2)}{\log 2}}.
	  \end{align*}This completes the proof. 
	  \end{proof}
	 \end{lem}

		\begin{cor}\label{aaa}
	 Suppose that $(Y,\nu_Y)\in D_{C,R}$ for some $R>0$ and $C:(0,R]\to [1,+\infty)$. Then $B_Y(x,r)$ is compact for any $x\in Y$
	 and $r>0$ with $2r\leq R$.
	 \begin{proof}The proof is by contradiction. Suppose that $B_Y(x,r)$ is not compact. Then, there exist $\varepsilon >0$ with
	  $\varepsilon \leq \min\{R-r,r\}$ and infinite 3$\varepsilon$-separated set 
	  $\{ x_i  \}_{i=1}^{\infty} \subseteq B_Y(x,r)$. By using Lemma \ref{aaaa}, we have
	  \begin{align*}
	   +\infty > \nu_Y\big(  B_Y (x,R)          \big)\geq \sum_{i=1}^{\infty}\nu_Y \big(B_Y  ( x_i, \varepsilon  )    \big)\geq
	   \sum_{i=1}^{\infty} \frac{1}{C(\varepsilon,r)^2}\Big( \frac{\varepsilon}{r} \Big)^{C(\varepsilon,r)}\nu_Y\big(    B_Y
	   (x,r)     \big) = +\infty,
	   \end{align*}which implies a contradiction. This completes the proof.
	  \end{proof}
		\end{cor}

	\subsection{Separation and concentration}
	In this subsection, we prove several results in \cite{gromov} because we find no proof anywhere.

	Let $(X,\dist)$ be a metric space. For $x\in X, r>0$, and $A,B\subseteq X$, we put
	\begin{align*}
	 \dist(A,B):=\inf \{ \dist(a,b) \mid a\in A,b\in B      \},\ \dist(x,A):=\dist(\{ x  \},A),\ A_r:=\{ x\in X \mid
	 \dist(x,A)\leq r    \}.
	 \end{align*}
		\begin{dfn}\upshape
	Let $(X, \dist , \mu )$ be an mm-space. For any $\kappa_0, \kappa_1, \cdots
	 ,\kappa_N \in \mathbb{R}$, we define 
	 \begin{align*}
	  \sep (X ; \kappa_0 , \cdots ,\kappa_N)=\ & \sep (\mu ; \kappa_0, \cdots
	  , \kappa_N)\\
	 : =\ & \sup \{  \min_{i\neq j} \dist(X_i, X_j)   \mid   X_0 , \cdots , X_N
	   \text{ {\rm are Borel subsets of }} X \\ & \hspace{4.9cm}\text{{\rm such that }}
	  \mu(X_i) \geq \kappa_i \text{ {\rm for any }}i \}, 
	 \end{align*}and call it the \emph{separation distance} of $X$.
		\end{dfn}

	The proof of the following lemma is easy and we omit the proof.
	\begin{lem}[cf. \cite{gromov}]\label{neko}
	 Let $(X, \dist_X, \mu_X)$ and $(Y, \dist_Y , \mu_Y) $ be two mm-spaces. Suppose that an $1$-Lipschitz map $f : X\to Y$
	 satisfies $f_{\ast}(\mu_X)= \mu_Y$. Then we have 
	 \begin{align*}
	  \sep (Y; \kappa_0, \cdots , \kappa_N) \leq \sep (X ; \kappa_0, \cdots , \kappa_N).
	 \end{align*}
	\end{lem}

	\begin{lem}[cf. \cite{gromov}]\label{noranoraneko}Let $(X,\dist,\mu)$ be an mm-space and $\kappa,\kappa' >0$ with $\kappa > \kappa'$. Then we have
	 \begin{align*}
	  \diam (X\oblip \mathbb{R},m-\kappa')\geq \sep (X;\kappa,\kappa).
	  \end{align*}
	 \begin{proof}Let $X_1,X_2\subseteq
	  X$ be two closed subsets such that $\mu (X_1) ,  \mu(X_2)\geq \kappa$. We define a function $f:X\to \mathbb{R}$ by $f(x):=
	  \dist(x,X_1)$. Let us show $\diam (f_{\ast}(\mu),m-\kappa')\geq  \dist(X_1,X_2)$ by contradiction. Suppose $\diam
	  (f_{\ast}(\mu),m-\kappa')< \dist(X_1,X_2)$. There exists a closed subset $X_0
	  \subseteq [0,+\infty)$ such that $\diam X_0 < \dist (X_1,X_2)$ and $\mu(f^{-1}(X_0))\geq m-\kappa'$. If 
	  $f^{-1}(X_0)\cap X_1 = \emptyset $, we have a contradiction since
	  \begin{align*}
	   \mu(f^{-1}(X_0)\cup X_1 )=\mu(f^{-1}(X_0))+ \mu(X_1) \geq (m-\kappa')+\kappa>m.
	   \end{align*}
	  In the same way we have $f^{-1}(X_0) \cap X_2 \neq
	  \emptyset$. Take a point $x_1 \in f^{-1}(X_0)\cap X_1$. $f(x_1)=\dist(x_1,X_1)=0\in X_0$ implies that $X_0 \subseteq
	  [0,\diam X_0]$. Therefore, we have $f^{-1}(X_0)\subseteq (X_1)_{\diam X_0}$, which yields $f^{-1}(X_0)\cap
	  X_2=\emptyset$ since $\diam X_0 < \dist (X_1,X_2)$. This is a contradiction since $f^{-1}(X_0) \cap X_2 \neq \emptyset$. As
	  a consequence, we obtain $\diam (f_{\ast}(\mu),m-\kappa')\geq \dist(X_1,X_2)$ and this completes the proof of the lemma. 
	  \end{proof}
	 \end{lem}

	\begin{rem}\upshape In \cite{gromov}, Lemma \ref{noranoraneko} is stated as $\kappa =\kappa'$, but that is not true in general. For example, let
	 $X:=\{ x_1 , x_2\}$, $\dist (x_1,x_2):=1$, and $\mu (\{ x_1\})=\mu (\{  x_2   \}):= 1/2$. Putting $\kappa =\kappa'=1/2$, we have
	 $\diam (X\oblip \mathbb{R},1/2)=0$ and $\sep
	 (X;1/2,1/2)=1$.
	 \end{rem}
	We denote by $\supp \mu$ the support of a Borel measure $\mu$. 
	\begin{lem}[cf. \cite{gromov}]Suppose that $\supp \mu$ is connected. Then, for any $\kappa >0$ we have
	 \begin{align*}
	  \diam (X \oblip \mathbb{R},m-\kappa)\geq \sep (X; \kappa, \kappa).
	  \end{align*}
	 \begin{proof}
	  Let $X_1,X_2\subseteq
	  X$ be two closed subsets such that $\mu (X_1)\geq \kappa$ and $\mu(X_2)\geq \kappa$. Define a function $f:X\to \mathbb{R}$
	  by $f(x):=
	  \dist(x,X_1)$. We will show $\diam (f_{\ast}(\mu),m-\kappa)\geq  \dist(X_1,X_2)$ by contradiction. Supposing that $\diam
	  (f_{\ast}(\mu),m-\kappa)< \dist(X_1,X_2)$, there exists a closed subset $X_0
	  \subseteq [0,+\infty)$ such that $\diam X_0 < \dist (X_1,X_2)$ and $\mu(f^{-1}(X_0))\geq m-\kappa$. If
	  $f^{-1}(X_0)\cap X_1 = \emptyset $, we have
	  \begin{align*}
	   \mu(f^{-1}(X_0)\cup X_1 )=\mu(f^{-1}(X_0))+ \mu(X_1) \geq (m-\kappa)+\kappa = m,
	   \end{align*}which implies $\supp \mu \subseteq f^{-1}(X_0)\cup X_1$. This is a contradiction since $\supp \mu$ is
	  connected. 
	  In the same way, we have $f^{-1}(X_0) \cap X_2 \neq
	  \emptyset$. Picking $x_1 \in f^{-1}(X_0)\cap X_1$, we get $f(x_1)=\dist(x_1,X_1)=0\in X_0$, which yields $X_0 \subseteq
	  [0,\diam X_0]$. Hence, we have $f^{-1}(X_0)\subseteq (X_1)_{\diam X_0}$, which implies $f^{-1}(X_0)\cap
	  X_2=\emptyset$ since 
	  $(X_1)_{\diam X_0}\cap X_2= \emptyset$. This is a contradiction because $f^{-1}(X_0) \cap X_2 \neq \emptyset$. As a result,
	  we obtain $\diam (f_{\ast}(\mu),m-\kappa)\geq \dist(X_1,X_2)$, which completes the proof of the lemma. 
	  \end{proof}
	 \end{lem}
	\begin{lem}[cf. \cite{gromov}]\label{sake}Let $\nu$ be a Borel measure on $\mathbb{R}$ with $m:= \nu
	 (\mathbb{R})<+\infty$. Then, for any $\kappa >0$ we have
	 \begin{align*}
	  \diam (\nu, m-2\kappa)\leq \sep (\nu; \kappa, \kappa).
	  \end{align*}
	 \begin{proof}Put $a_0:= \sup \{ a\in \mathbb{R} \mid \nu \big( (-\infty,a) \big)\leq \kappa \}$ and $b_0:= \inf \{ b\in \mathbb{R}
	  \mid \nu \big( (b, + \infty)  \big)\leq \kappa\}$. Then, we have $a_0\leq b_0$ and 
	  \begin{align*}
	   \kappa \leq \ &\lim_{\varepsilon  \downarrow 0} \nu \big((-\infty, a_0+\varepsilon)\big)= \nu \big((-\infty, a_0]\big), \\
	   \kappa \leq \ & \lim_{\varepsilon  \downarrow 0} \nu \big((b_0-\varepsilon,+\infty)\big)= \nu \big([b_0,+ \infty)\big).
	   \end{align*}$\nu \big( (-\infty, a_0) \big)\leq \kappa$ and $\nu \big( (b_0, +\infty) \big)\leq
	  \kappa$ imply $\nu ([a_0,b_0])\geq m-2\kappa$. Therefore, indicating by $\dist_{\mathbb{R}}$ the usual Euclidean distance, we obtain
	  \begin{align*}
	   \diam (\nu, m-2\kappa)\leq \diam ([a_0,b_0])=b_0-a_0 =\dist_{\mathbb{R}} \big((-\infty,a_0],[b_0,+\infty) \big)\leq \sep (\nu;\kappa,\kappa).
	   \end{align*}This completes the proof.
	  \end{proof}
	 \end{lem}
	\begin{cor}[cf. \cite{gromov}]\label{noraneko}
	 For any $\kappa >0$, we have
	 \begin{align*}
	  \sep (X;\kappa,\kappa)\geq \diam (X \oblip \mathbb{R},m-2\kappa).
	  \end{align*}
	 \begin{proof}Let $f:X\to \mathbb{R}$ be an arbitrary $1$-Lipschitz function. From Lemma $\ref{neko}$ and Lemma $\ref{sake}$,
	  we have $\sep (X;\kappa,\kappa)\geq \sep (f_{\ast}(\mu);\kappa,\kappa) \geq \diam (f_{\ast}(\mu),m-2\kappa)$. This completes the proof. 
	  \end{proof}
	 \end{cor}
	Combining Lemma \ref{noranoraneko} and Corollary \ref{noraneko} we obtain the following corollary. 
	\begin{cor}[cf. \cite{gromov}]A sequence $\{ X_n \}_{n=1}^{\infty}$ of mm-spaces is a L\'{e}vy family if and only if
	 $\sep (X_n; \kappa, \kappa)\to 0$ as $n\to \infty$ for any $\kappa >0$. 
	 \end{cor}

	\section{Proof of the Main Theorem}

	\begin{proof}[Proof of Theorem $\ref{nobunaga}$]
	Let $ \{  (Y_n,\nu_{Y_n})     \}_{n=1}^{\infty}$ be any sequence of $D_{C,R}$ and $\{ f_n :X_n \to Y_n     \}_{n=1}^{\infty}$
	 any sequence of $1$-Lipschitz maps. Given any $\varepsilon > 0$ with $32\varepsilon  \leq 3R$, it suffies to show that $\diam (f_{n
	 \ast}(\mu_n),m_n-\kappa)\leq 6\varepsilon$ for any $n$ by choosing a subsequence. The claim obviously
	 holds in the case of $\limsup\limits_{n\to \infty} m_n =0$, so we assume that $\inf\limits_{n\in \mathbb{N}} m_n
	 >0$. Take a maximal $\varepsilon$-separated set $\{ \xi_{\alpha}^n    \}_{\alpha \in \mathcal{A}_n}\subseteq Y_n$ for each
	 $n\in \mathbb{N}$.
	 \begin{claim}\label{mihimaru}For any $n\in \mathbb{N}$ and $\alpha \in \mathcal{A}_n$, we have
	  \begin{align*}
	   \card  ( \{ \beta \in \mathcal{A}_n \mid \xi_{\beta}^n \in B_{Y_n}(\xi_{\alpha}^n,5\varepsilon)       \})\leq
	   2^{4C(\varepsilon/3, 16\varepsilon/3)} C\Big(\frac{\varepsilon}{3}, \frac{16\varepsilon}{3}\Big)^2.
	   \end{align*}
	  \begin{proof}By Cororally \ref{aaa}, the set $\{ \beta \in \mathcal{A}_n \mid \xi_{\beta}^n \in
	   B_{Y_n}(\xi_{\alpha}^n,5\varepsilon)       \}$ is finite. Let $\{ {\beta}_1 , {\beta}_2 , \cdots
	   {\beta}_k      \}:= \{ \beta \in \mathcal{A}_n \mid \xi_{\beta}^n \in B_{Y_n}(\xi_{\alpha}^n,5\varepsilon)       \}$ and
	   take $j\in \{  1,2,\cdots ,k   \}$ such that $\nu_{Y_n}\big(
	     B_{Y_n}(\xi_{\beta_{j}}^n,\varepsilon/3) \big) =  \min \{    \nu_{Y_n}\big(B_{Y_n}(\xi_{\beta_{\ell}}^n,\varepsilon/3)\big) \mid \ell = 1,2, \cdots ,k
	   \}$. Since 
	   \begin{align*}
		\nu_{Y_n}\Big(    B_{Y_n}\Big(\xi_{\alpha}^n ,\frac{16\varepsilon}{3}\Big)       \Big) \geq   \sum_{\ell=1}^k \nu_{Y_n} \Big(
		B_{Y_n}\Big(\xi_{{\beta}_{\ell}}^n,\frac{\varepsilon}{3}\Big) \Big) \geq k \nu_{Y_n}\Big(B_{Y_n}\Big(\xi_{\beta_j}^n,\frac{\varepsilon}{3}\Big)\Big),
		\end{align*}combining this and Lemma \ref{aaaa}, we have
	   \begin{align*}
		k\leq  \frac{\nu_{Y_n}\big(   B_{Y_n}(\xi_{\alpha}^n ,16\varepsilon/3)             \big)}{\nu_{Y_n} \big(   B_{Y_n}(\xi_{\beta_j}^n,\varepsilon /3)
		\big)}\leq \Big(\frac{16\varepsilon/3}{\varepsilon /3}\Big)^{C(\varepsilon /3
		,16\varepsilon /3)}C\Big(\frac{\varepsilon}{3},\frac{16\varepsilon}{3}\Big)^2= 2^{4C(\varepsilon/{3}, 16\varepsilon/3)}C\Big(\frac{\varepsilon}{3}, \frac{16\varepsilon}{3}\Big)^2.
		\end{align*}This completes the proof of Claim \ref{mihimaru}.
	   \end{proof}
	  \end{claim}
	 By Claim \ref{mihimaru}, for each $n\in \mathbb{N}$ there exists ${\alpha}_n\in \mathcal{A}_n$ such that
	 \begin{align*}
	  k_n:=\card (\{  \beta \in \mathcal{A}_n\mid \xi_{\beta}^n  \in B_{Y_n}(\xi_{{\alpha}_n}^n,5\varepsilon)     \})=
	 \sup\limits_{\alpha \in \mathcal{A}_n} \card (\{   \beta \in \mathcal{A}_n \mid \xi_{\beta}^n \in B_{Y_n}(\xi_{\alpha}^n,5\varepsilon)              \}).
	  \end{align*}By taking a subsequence, we get $k_n \equiv k$ for any $n$. Put $\{ {\beta}_1^n,{\beta}_2^n, \cdots
	 ,{\beta}_k^n     \}:= \{  \beta \in \mathcal{A}_n \mid
	 \xi_{\beta}^n \in B(\xi_{{\alpha}_n}^n,5\varepsilon)   \}$. We take $J_1^n \subseteq \{  \xi_{\alpha}^n   \}_{\alpha \in
	 \mathcal{A}_n}$ which is maximal with respect to the properties that $J_1^n$ is
	 $5\varepsilon$-separated and $\xi_{{\beta_1^n}}^n \in J_1^n$, $\xi_{{\beta}_2^n}^n \notin J_1^n$, $\cdots$ ,
	 $\xi_{{\beta}_k^n}^n\notin J_1^n$. Next, we take $J_2^n
	 \subseteq \{   \xi_{\alpha}^n   \}_{\alpha \in \mathcal{A}_n} \setminus J_1^n$ which is maximal with respect to the properties that $J_2^n$ is $5\varepsilon$-separated and $\xi_{{\beta}_2^n}^n \in J_2^n$, $\xi_{{\beta}_3^n}^n \notin J_2^n$,
	 $\cdots$, $\xi_{{\beta}_k^n}^n \notin J_2^n$. In the same way, we pick $J_3^n \subseteq \{  \xi_{\alpha}^n
	 \}_{\alpha \in \mathcal{A}_n}\setminus (J_1^n \cup J_2^n)$, $\cdots$, $J_{\alpha}^n \subseteq \{  \xi_{\alpha}^n
	 \}_{\alpha \in \mathcal{A}_n} \setminus
	 (J_1^n \cup J_2^n \cup \cdots \cup J_{k-1}^n)$. Then we have
	 \begin{claim}\label{onemu}$\{  \xi_{\alpha}^n    \}_{\alpha \in \mathcal{A}_n}= J_1^n \cup J_2^n \cup \cdots \cup J_k^n$ for each $n\in \mathbb{N}$.
	  \begin{proof}The proof is by contradiction. Let us suppose that $\xi_{\alpha}^n \notin J_1^n \cup J_2^n \cup \cdots \cup J_k^n$. 
	   Since $J_i^n$ is maximal  for each $i=1,2, \cdots ,k$, there exists 
	   $\xi_{{\gamma}_i}^n \in J_i^n$ such that $\dist_{Y_n} (\xi_{\alpha}^n , \xi_{{\gamma}_i}^n)< 5\varepsilon$ and $\xi_{\alpha}^n \neq
	   \xi_{{\gamma}_i}^n$.
	   Therefore, we have
	   \begin{align*}
		k+1\leq \card(\{ \xi_{\alpha}^n,\xi_{{\gamma}_1}^n ,\xi_{{\gamma}_2}^n , \cdots ,\xi_{{\gamma}_k}^n  \})\leq \card (\{
		\beta \in \mathcal{A}_n \mid \xi_{\beta}^n \in
		B_{Y_n}(\xi_{\alpha}^n,5\varepsilon)            \})\leq k,
		\end{align*}which is a contradiction. 
	   \end{proof}
	  \end{claim}
	 By Claim \ref{onemu}, we have $Y_n = \bigcup\limits_{i=1}^k \bigcup\limits_{\xi_{\alpha}^n \in J_i^n} B_{Y_n}(\xi_{\alpha}^n ,\varepsilon)
	 $. Therefore, by taking a subsequence, there exists $j_0\in \mathbb{N}$ such that $1\leq j_0 \leq k$ and 
	 \begin{align*}
	  f_{n \ast}(\mu_n)\Big( \bigcup_{\xi_{\alpha}^n\in J_{j_0}^n} B_{Y_n}(\xi_{\alpha}^n,\varepsilon)   \Big)\geq \frac{m_n}{k}
	  \end{align*}for any $n\in \mathbb{N}$. Then we obtain
	 \begin{claim}\label{cla}
	  \begin{align*}
	   f_{n \ast}(\mu_n)\Big( Y_n \setminus \bigcup_{\xi_{\alpha}^n\in J_{j_0}^n} B_{Y_n}(\xi_{\alpha}^n,2\varepsilon)   \Big) \to 0
	   \end{align*}as $n\to \infty$.
	  \begin{proof}If $f_{n \ast}(\mu_n)\Big( Y_n \setminus \bigcup\limits_{{\xi}_{\alpha}^n\in J_{j_0}^n} B_{Y_n}(\xi_{\alpha}^n,2\varepsilon)
	   \Big) \geq C'$ for a constant $C'>0$ and infinitely many $n\in \mathbb{N}$, then we have
	   \begin{align*}
		\varepsilon \leq \dist_{Y_n}\Big( Y_n \setminus \bigcup\limits_{{\xi}_{\alpha}^n \in J_{j_0}^n} B_{Y_n}(\xi_{\alpha}^n,2\varepsilon),
		\bigcup\limits_{{\xi}_{\alpha}^n \in J_{j_0}^n} B_{Y_n}(\xi_{\alpha}^n,\varepsilon)  \Big)\leq \ &\sep \Big( f_{n \ast}(\mu_n); C',\frac{1}{k}\inf_{n\in \mathbb{N}}
		m_n          \Big)\\
		\leq \ &\sep \Big( \mu_n; C',\frac{1}{k}\inf_{n\in \mathbb{N}} m_n          \Big).
		\end{align*}This is a contradiction, since the right-hand side of the above inequality converges to $0$ as $n\to \infty$. We
	   have the claim.
	   \end{proof}
	  \end{claim}
	 \begin{claim}For any suffiecient large $n\in \mathbb{N}$ there exists $\xi_{{\gamma}_n}^n \in J_{j_0}^n$ such that
	  \begin{align*}
	   f_{n \ast}(\mu_n)\big(  B_{Y_n}(\xi_{{\gamma}_n}^n,2\varepsilon)           \big)\geq \frac{1}{6} \inf_{n\in \mathbb{N}}m_n.
	   \end{align*}
	  \begin{proof}Let us prove the claim by contradiction. Suppose that
	   \begin{align*}
		f_{n \ast}(\mu_n)\big(  B_{Y_n}(\xi_{\alpha}^n,2\varepsilon)           \big)< \frac{1}{6} \inf_{n\in \mathbb{N}}m_n
		\end{align*}for infinitely many $n\in \mathbb{N}$ and any $\xi_{\alpha}^n \in J_{j_0}^n$. By Claim \ref{cla}, there
	   exist $n_0\in \mathbb{N}$ such that
	   \begin{align*}
		f_{n \ast}(\mu_n) \Big( \bigcup_{\xi_{\alpha}^n\in J_{j_0}^n} B_{Y_n}(\xi_{\alpha}^n,2\varepsilon) \Big)\geq
		\frac{5}{6}\inf_{n\in \mathbb{N}} m_n.
		\end{align*}for any $n\in \mathbb{N}$ with $n\geq n_0$. From the assumption, if $n\geq n_0$ we have $J_n' \subseteq J_n$ such that
	   \begin{align*}
		 \frac{1}{6} \inf_{n\in \mathbb{N}} m_n \leq f_{n \ast}(\mu_n)\Big(\bigcup_{\xi_{\alpha}^n \in
		J_n'}B_{Y_n}(\xi_{\alpha}^n,2\varepsilon)\Big)\leq \frac{1}{3} \inf_{n\in \mathbb{N}}m_n.
		\end{align*}Hence, by putting $J_n'':= J_n \setminus J_n'$ we have
	   \begin{align*}
		\varepsilon \leq \dist_{Y_n}\Big(    \bigcup_{\xi_{\alpha}^n \in J_n'} B_{Y_n}(\xi_{\alpha}^n,2\varepsilon) ,
		\bigcup_{\xi_{\alpha}^n \in J_n''} B_{Y_n}(\xi_{\alpha}^n ,2\varepsilon)         \Big)\leq \ &\sep \Big(  f_{n \ast}(\mu_n);
		\frac{1}{6} \inf_{n\in \mathbb{N}} m_n , \frac{1}{2} \inf_{n\in \mathbb{N}}  m_n  \Big)\\\leq \ &\sep \Big( \mu_n ; \frac{1}{6} \inf_{n\in \mathbb{N}} m_n , \frac{1}{2} \inf_{n\in \mathbb{N}}  m_n     \Big),
		\end{align*}which is a contradiction since the right-hand side of the above inequality converges to $0$ as $n\to \infty$. 
	   \end{proof}
	  \end{claim}
	 \begin{claim}\label{nemuinoda}We have
	  \begin{align*}
	   f_{n \ast}(\mu_n) \big( Y_n \setminus B_{Y_n}(\xi_{{\gamma}_n}^n,3\varepsilon)           \big)\to 0
	   \end{align*}as $n\to \infty$. 
	  \begin{proof}The claim immediately follows from the same proof of Claim \ref{cla}. 
	   \end{proof}
	  \end{claim}
	 By Claim \ref{nemuinoda}, for any suffiecient large $n\in \mathbb{N}$ we have
	 \begin{align*}
	  f_{n \ast}(\mu_n) \big( Y_n \setminus B_{Y_n}(\xi_{{\gamma}_n}^n,3\varepsilon)           \big)\leq \kappa,
	  \end{align*}which implies $\diam (f_{n \ast}(\mu_n),m_n -\kappa)\leq \diam
	 B_{Y_n}(\xi_{{\gamma}_n}^n,3\varepsilon)\leq 6\varepsilon$. This completes the proof of Theorem \ref{nobunaga}.
	\end{proof}

	\begin{ack}\upshape
	 The author would like to thank Professor Takashi Shioya for valuable discussions. 
	 \end{ack}


\begin{thebibliography}{3}
	 \bibitem{ambro}L. Ambrosio and P. Tilli,  {\it Topics on analysis in metric spaces}, Oxford Lecture Series in Mathematics and its Applications, 25. Oxford University Press, Oxford, 2004. 
	\bibitem{funano}K. Funano, {\it Observable concentration of mm-spaces with hyperbolic spaces as screens}, preprint, 2007.
	 \bibitem{gromov}M. Gromov,  {\it Metric structures for Riemannian and non-Riemannian spaces}, Based on the 1981 French
			 original, With appendices by M. Katz, P. Pansu and S. Semmes. Translated from the French by Sean Michael
			 Bates. Progress in Mathematics, 152. Birkh\"{a}user Boston, Inc., Boston, MA, 1999.
	\bibitem{ledoux}M. Ledoux, {\it The concentration of measure phenomenon}, Mathematical Surveys and Monographs, 89. American
			 Mathematical Society, Providence, RI, 2001.
	\end{thebibliography}
	\end{document}